\documentclass[a4paper,12pt,reqno]{amsart}

\usepackage{amsmath}
\usepackage{amsthm}
\usepackage{amssymb}
\usepackage{latexsym}
\usepackage{setspace}

\newtheorem{thm}{Theorem}[section]

\newtheorem{prop}{Proposition}[section]
\newtheorem{lem}{Lemma}[section]

\theoremstyle{definition}

\newtheorem{rem}{Remark}[section]

\numberwithin{equation}{section}
\allowdisplaybreaks[4]

\newcommand{\pt}{\partial}
\newcommand{\rre}{{\mathbb R}}

\title[nonlinear Schr\"odinger equation with 
delta potential]{
Final state problem for the cubic nonlinear\\
Schr\"odinger equation with repulsive \\
delta potential} 
\author[Jun-ichi Segata]{Jun-ichi Segata}



\date{}

\begin{document}

\maketitle

\centerline{Mathematical Institute, Tohoku University,}
\centerline{Aoba, Sendai 980-8578, Japan}
\centerline{segata@math.tohoku.ac.jp}

\vskip5mm
\noindent
{\bf Abstract.} 
We consider the asymptotic behavior in time of solutions to 
the cubic nonlinear Schr\"odinger equation with repulsive 
delta potential 
($\delta$-NLS). 
We shall prove that for 
a given asymptotic profile $u_{ap}$, there exists a solution 
$u$ to ($\delta$-NLS) which converges to $u_{ap}$ in $L^{2}(\rre)$ 
as $t\to\infty$. 
To show this result we exploit the distorted Fourier transform 
associated to the Schr\"odinger equation with delta potential.

\section{Introduction} \label{sec:intro}

We consider the nonlinear 
Schr\"{o}dinger equation 
with delta potential: 
\begin{eqnarray}
i\pt_tu+\frac12\pt_x^2u
-q\delta u
=\lambda|u|^2u,\qquad t,x\in\rre,
\label{NLS}
\end{eqnarray}
where  $q\neq0$ and $\lambda\neq0$ are real constants, and 
$\delta$ is the Dirac measure at the origin. 

The cubic nonlinear Schr\"odinger equation with localized 
potential provides a simpler model describing the resonant 
nonlinear propagation of light through optical wave guides 
with localized defects, see Goodman-Holmes-Weinstein \cite{GHW}. 
The delta potential is called ``repulsive'' for $q>0$ 
and ``attractive'' for $q<0$. 

Along with the 
nonlinear Schr\"odinger equation, 
(\ref{NLS}) has a standing wave solution. 
The stability of the soliton 
is studied in Goodman-Holmes-Weinstein \cite{GHW}, 
Fukuizumi-Ohta-Ozawa \cite{FOO} by the variational arguments, 
and Deift-Park \cite{DP} via the nonlinear steepest-descent 
method for Riemann-Hilbert problems. 

Eq. (\ref{NLS}) is currently being intensively investigated 
in the point of view of 
the interaction between the soliton and the delta potential 
\cite{DH,GHW,HMZ1,HMZ2,HZ}. 
Holmer-Zworski \cite{HZ} studied the behavior of slow solitons 
under the flow by (\ref{NLS}). Later on, 
Holmer-Marzuola-Zworski \cite{HMZ1,HMZ2} and 
Datchev-Holmer \cite{DH} have studied  
the splitting of the fast solitons by the delta potential. 
It was shown that the high velocity incoming soliton is 
split into a transmitted component and a reflected component 
for both repulsive \cite{HMZ2} and attractive \cite{DH} cases, 
see \cite{GHW,HMZ1} for the numerical results. 

In this paper we study the long time behavior of solutions 
to (\ref{NLS}). More precisely, we shall prove that for 
a given asymptotic profile $u_{ap}$, there exists a solution 
$u$ to (\ref{NLS}) which converges to $u_{ap}$ in $L^{2}(\rre)$ 
as $t\to\infty$. So we consider the final state problem (\ref{NLS}) and 
\begin{eqnarray}
\lim_{t\to+\infty}(u(t)-u_{ap})=0\qquad in\ L^{2}.
\label{ac}
\end{eqnarray}
Although we consider the behavior of solution to 
(\ref{NLS}) as $t\to\infty$, we can treat behavior of solution 
as $t\to-\infty$ in a similar fashion. 

Before we state our result, 
we summarize known results for the scattering and 
asymptotic behavior of solution to the  
nonlinear Schr\"odinger equation. 
There is a large literatures on the 
scattering theory of the solution to 
nonlinear Schr\"odinger equation 
starting with the pioneering work of 
Lin-Strauss \cite{LS} and Ginibre-Velo \cite{GV3}. 
We are concerned here only with 
the one dimensional case. 

Let us consider the nonlinear Schr\"odinger equation with 
a potential
\begin{eqnarray}
i\pt_tv+\frac12\pt_x^2v-V(x)v=\mu|v|^{p-1}v,\qquad t,x\in\rre,
\label{aNLS}
\end{eqnarray}
where $p>1$, $\mu\neq0$ and $V$ is a real valued function. 
For the case $V\equiv0$, it is known 
that the power $p=3$ is the borderline between the short and 
long range scattering theory, namely, 
the solution to (\ref{aNLS}) scatters to the solution 
to the free Schr\"odinger equation for $p>3$    
and no solution of (\ref{aNLS}) has a scattering 
state for $p\le3$, see Barab \cite{B} and Tsutsumi-Yajima  
\cite{TY}. Therefore the solution of the 
nonlinear equation (\ref{NLS}) with $p\le3$ may be have 
a different asymptotic profile from the solution of the 
free Schr\"odinger equation. 

In this regard, Ozawa \cite{O}
proved that for a given small final data $\phi_{+}$, 
there exists a solution $v$ to 
(\ref{aNLS}) with $p=3$ satisfying
\begin{eqnarray}
v(t)\sim v_{ap}(t)=t^{-1/2}{{\mathcal F}}[\phi_{+}](\frac{x}{t})
\exp\left(\frac{ix^{2}}{2t}-i\lambda|{{\mathcal F}}[\phi_{+}]
(\frac{x}{t})|^{2}\log t
-i\frac{\pi}{4}\right)\label{D}
\end{eqnarray}
as $t\to+\infty$, where ${{\mathcal F}}$ is the 
Fourier transform. 
Note that $v_{ap}$ is the leading term of $e^{(i/2)t\pt_{x}^{2}}
\phi_{+}$ with the phase modification. 
Later on, Hayashi-Naumkin \cite{HN1} proved that 
for a given small initial data $v_{0}$, 
there exists a unique solution 
$v$ to (\ref{aNLS}) with $p=3$ 
satisfying $v(0,x)=v_{0}(x)$ and (\ref{D}) for some $\phi_{+}$. 
See also Hayashi-Naumkin \cite{HN2} for the improvement 
of those results. 

For the case where $V\not\equiv0$, 
Cuccagna-Georgiev-Visciglia \cite{CGV} 
have recently shown that the solution to (\ref{aNLS}) 
scatters to the solution of the free Schr\"odinger equation
under the suitable assumptions on $V$ and $p>3$. 
As far as we know, there is no result on 
the long range scattering for (\ref{aNLS}) with 
$V\not\equiv0$. 

We briefly explain why the 
critical exponent $p$ is three for the one 
dimensional Schr\"odinger equation. 
Roughly speaking, the 
nonlinearity of (\ref{aNLS}) 
is short range if 
$L^2$ norm of the nonlinear 
term is integrable on time interval $[1,\infty)$. 
Since the pointwise decay of the solution to the 
one dimensional free Schr\"odinger equation is 
${\mathcal O}(t^{-1/2})$ as $t\to\infty$, 
the $L^2$ norm of the nonlinear term $|v|^{p-1}v$ decays like 
${\mathcal O}(t^{-(p-1)/2})$. Since the integral 
$\int_1^{\infty}t^{-(p-1)/2}dt$ is finite 
if and only if $p>3$, the exponent $p=3$ appears 
as the borderline between the short and 
long range scattering theory. 
For further results
on the scattering theory on the nonlinear 
Schr\"odinger equation, see e.g., 
Cazenave \cite[Chapter 7]{C}.

Let us return to the equation (\ref{NLS}). As with 
(\ref{aNLS}), the pointwise decay of the solution to the 
free Schr\"odinger equation with delta potential is ${\mathcal O}
(t^{-1/2})$, see proof of Lemma \ref{S} below.  
Hence we expect that (\ref{NLS}) fall into the long range 
case. In fact, 
we obtain the following results for the final state problem 
(\ref{NLS})-(\ref{ac}).

\begin{thm}\label{main}
Let $q\ge0$. There exists $\varepsilon>0$ and $C>0$ 
with the following properties: 
For any $\phi_+\in{{\mathcal S}}'(\rre)$  
with $\|(1+|x|)\phi_+\|_{L_{x}^{2}}<\varepsilon$ 
there exists a unique global solution $u\in C(\rre;L_x^2(\rre))
\cap L_{loc}^4(\rre;L_x^{\infty}(\rre))$ 
to (\ref{NLS}) satisfying
\begin{eqnarray*}
\sup_{t\ge1}t^{\alpha}\left\{
\|u(t)-u_{ap}(t)\|_{L_x^2}
+\left(\int_{t}^{\infty}\|u(\tau)-u_{ap}(\tau)
\|_{L_x^{\infty}}^{4}d\tau\right)^{1/4}\right\}\le C,
\end{eqnarray*}
where $1/4<\alpha<1/2$ and $u_{ap}$ is given by 
\begin{eqnarray*}
u_{ap}(t,x)=t^{-1/2}{{\mathcal F}}_{q}[\phi_{+}](\frac{x}{t})
\exp\left(\frac{ix^{2}}{2t}-i\lambda|{{\mathcal F}}_{q}[\phi_{+}]
(\frac{x}{t})|^{2}\log t
-i\frac{\pi}{4}\right)
\end{eqnarray*}
and ${{\mathcal F}}_{q}$ is a 
distorted Fourier transform associated to 
Schr\"{o}dinger operator with delta potential: 
\begin{eqnarray*}
{{\mathcal F}}_{q}[\phi](\xi)
&=&\left\{
\begin{array}{l}
\medskip\displaystyle{
{{\mathcal F}}[\phi(x)-q1_{-}(x)e^{qx}
\int_{x}^{-x}e^{q|y|}\phi(y)dy](\xi)
\ \qquad\text{if}\ \ \xi\ge0,}
\nonumber\\
\displaystyle{
{{\mathcal F}}[\phi(x)-q1_{+}(x)e^{-qx}
\int_{-x}^{x}e^{q|y|}\phi(y)dy](\xi)
\qquad\text{if}\ \ \xi<0.}
\end{array}
\right.
\end{eqnarray*}
\end{thm}

\vskip3mm

\begin{rem} 
In this paper we consider the repulsive case only. 
The reason is due to a difference in the spectrum properties 
of the Hamiltonian $H_{q}=-(1/2)\pt_{x}^{2}+q\delta(x)$.  
In fact, if $q\ge0$, then $H_{q}$ 
has no eigenvalues and if $q<0$, $H_{q}$ has precisely one negative, 
simple eigenvalue, see Section 2 below. 
Hence for the attractive case, 
we need to analyze $H_{q}$ by taking into 
account of the presence of the eigenfunction. 
\end{rem}

In proving Theorem \ref{main}, the one of key points is 
choice of the asymptotic profile $u_{ap}$. To 
explain how to choose $u_{ap}$, we employ the argument 
due to Ozawa \cite{O}. We first rewrite 
the final state problem (\ref{NLS})-(\ref{ac}) 
as the integral equation of Yang-Feldman type
\begin{eqnarray}
u(t)-u_{ap}(t)
=i\int_{t}^{\infty}
e^{-i(t-\tau)H_{q}}
\{\lambda|u|^{2}u-(i\pt_{t}-\frac12 H_{q})u\}(\tau)d\tau.
\label{INT}
\end{eqnarray}
We construct 
the solution to (\ref{INT}) via the contraction mapping 
principle in a suitable Banach space. To this end, we 
split the right hand side of (\ref{INT}) into the following 
two pieces
\begin{eqnarray*}
u(t)-u_{ap}(t)
&=&i\lambda\int_{t}^{\infty}
e^{-i(t-\tau)H_{q}}
\{|u|^{2}u-|u_{ap}|^{2}u_{ap}\}(\tau)d\tau\\
& &-i\int_{t}^{\infty}
e^{-i(t-\tau)H_{q}}
\{(i\pt_{t}-\frac12 H_{q})u_{ap}
-\lambda|u_{ap}|^{2}u_{ap}\}(\tau)d\tau.
\end{eqnarray*}
To apply the Banach fixed point theorem 
for our problem, we are required that the $L^{2}$ 
norm of the integrant appearing in the second term 
of the above equation decays faster than $t^{-5/4}$. 
To respond to the requirement, we take
\begin{eqnarray*}
u_{ap}(t,x)=t^{-1/2}{{\mathcal F}}_{q}[\phi_{+}](\frac{x}{t})
e^{\frac{ix^{2}}{2t}-\frac{i\pi}{4}+iS(t,\frac{x}{t})}.
\end{eqnarray*}
Note that $u_{ap}$ is the leading term of $e^{-itH_{q}}
\phi_{+}$ with the phase modification $S(t,x/t)$. 
Then an elementary calculation yields
\begin{eqnarray*}
(i\pt_{t}-\frac12 H_{q})u_{ap}
-\lambda|u_{ap}|^{2}u_{ap}
\sim (-\pt_{t}S(t,\frac{x}{t})
-\lambda t^{-1}|{{\mathcal F}}_{q}[\phi_{+}](\frac{x}{t})|^{2})
u_{ap}.
\end{eqnarray*}
Therefore taking $S(t,\xi)=-\lambda|{{\mathcal F}}_{q}[\phi_{+}](\xi)
|^{2}\log t$, 
we obtain the desired estimate 
and we can apply the fixed point argument. 
In this step we have to 
impose the strong assumption 
on the final data. To weaken the assumption 
on the final data, we use a modified version of 
the integral equation of Yang-Feldman type
which is introduced by Hayashi-Naumkin \cite{HN2}. 

The another crucial part of this paper is the 
derivation of the asymptotic formula:
\begin{eqnarray*}
(e^{-itH_{q}}\phi)(x)\sim t^{-1/2}{{\mathcal F}}_{q}[\phi](\frac{x}{t})
e^{\frac{ix^{2}}{2t}-\frac{i\pi}{4}}
\end{eqnarray*}
as $t\to+\infty$ in $L^{p}$ with $2\le p\le \infty$, 
see Proposition \ref{l3} below for the detail. 
To derive this formula, we use several nontrivial identities  
related to the distorted Fourier transform associated to $H_{q}$. 
The distorted Fourier transform and its inverse transform 
for $H_{q}$ are defined by
\begin{eqnarray*}
{{\mathcal F}}_{q}[\phi](\xi)
=\left\{
\begin{array}{l}
\medskip\displaystyle{
{{\mathcal F}}[{{\mathcal L}}_{+}[\phi](x)](\xi)
\ \quad\text{if}\ \xi\ge0,}\\
\displaystyle{
{{\mathcal F}}[{{\mathcal L}}_{-}[\phi](x)](\xi)
\quad\ \text{if}\ \xi<0,}
\end{array}
\right.
\quad
{{\mathcal F}}_{q}^{-1}[\phi](x)=
\left\{
\begin{array}{l}
\medskip\displaystyle{
{{\mathcal F}}_{q,+}^{-1}[\phi](x)
\ \quad\text{if}\ x\ge0,}\\
\displaystyle{
{{\mathcal F}}_{q,-}^{-1}[\phi](x)
\quad\ \text{if}\ x<0,}
\end{array}
\right.
\end{eqnarray*}
for some operators ${{\mathcal L}}_{\pm}$ and 
${{\mathcal F}}_{q,\pm}^{-1}$, see Section 3 for 
the explicit forms of those operators. 
Especially we derive the identities
\begin{eqnarray*}
{{\mathcal L}}_{+}
[{{\mathcal F}}_{q}^{-1}[\phi]](x)&=&
{{\mathcal F}}_{q,+}^{-1}[\phi](x),\\
{{\mathcal L}}_{-}
[{{\mathcal F}}_{q}^{-1}[\phi]](x)&=&
{{\mathcal F}}_{q,-}^{-1}[\phi](x)
\end{eqnarray*}
for any $x\in\rre$. Using the above identities and several 
properties of the distorted Fourier transform, 
we guarantee Proposition \ref{l3}. 

We introduce several 
notation and functional spaces. 
Let $1_{+}(x)$ and $1_{+}(x)$ be characteristic functions in $x$ 
on the intervals $[0,\infty)$ and $(-\infty,0)$. 
Let $\langle\xi\rangle=\sqrt{\xi^2+1}$ and 
$\langle i\pt_{x}\rangle=\sqrt{1-\pt_{x^{2}}}$.  
For $\alpha,\beta\in\rre$, we denote $H^{{\alpha,\beta}}$ 
the weighted Sobolev spaces
\begin{eqnarray*}
H^{\alpha,\beta}
&=&\{\phi\in{{\mathcal S}}'(\rre);\|\phi\|_{H^{\alpha,\beta}}<\infty\},\\
\|\phi\|_{H^{\alpha,\beta}}&=&
\|\langle x\rangle^{\beta}\langle i\pt_{x}\rangle^{\alpha}
\phi\|_{L_{x}^{2}}.
\end{eqnarray*}
For $1\le q,r\le\infty$, $L^q(t,\infty;L_x^r(\rre))$ is defined 
as follows:
\begin{eqnarray*}
L^q(t,\infty;L_x^r(\rre))&=&\{u\in{{\mathcal S}}'(\rre^2);
\|u\|_{L^q(t,\infty;L_x^r)}<\infty\},\\
\|u\|_{L^q(t,\infty;L_x^r)}&=&
\left(\int_t^{\infty}\|u(\tau)\|_{L_x^r}^qd\tau\right)^{1/q}.
\end{eqnarray*}
We use the notation $A\lesssim B$ to denote the estimate 
$A\le CB$ where $C$ is a positive constant.

In Section 2, we study several properties of 
the Hamiltonian $H_{q}$. 
Section 3 is devoted to introducing the distorted 
Fourier transform associated to $H_{q}$ and 
giving an asymptotic formula as $t\to\infty$ for 
the unitary group $e^{-itH_{q}}$. 
In Section 4, we conclude the proof of Theorem \ref{main} 
via the contraction mapping principle. 
In Appendix, we prove the key identities in Lemma \ref{ly}. 

\section{Linear estimates} \label{sec:linear}
 
We consider the Hamiltonian associated to 
linear Schr\"{o}dinger equation 
with repulsive delta potential: 
\begin{eqnarray*}
H_{q}=-\frac12\frac{d^{2}}{dx^{2}}+q\delta(x)\qquad\xi\in\rre.
\end{eqnarray*}
Let
\begin{eqnarray*}
\left\{
\begin{array}{l}
\medskip\displaystyle{D(H_{q})=
\{u\in H^{1}(\rre)\cap H^{2}(\rre\backslash\{0\});
u'(0+)-u'(0-)=2qu(0)\},}\\
\displaystyle{H_{q}u=
-\frac12\frac{d^{2}}{dx^{2}}.}
\end{array}
\right.
\end{eqnarray*}
Then $H_{q}$ is self-adjoint operator in $L^{2}(\rre)$, 
see \cite[Theorem 3.1.1 in Chapter 1.3]{AG}. 
Hence the Stone theorem yields that $H_{q}$ generates 
the $L^{2}$-unitary group $\{e^{-itH_{q}}\}_{t\in\rre}$. 
Furthermore, for $q\in\rre$, the essential spectrum 
of $H_{q}$ is purely absolutely continuous and 
\begin{eqnarray*}
\sigma_{ess}(H_{q})=\sigma_{ac}(H_{q})=[0,\infty),
\qquad \sigma_{sc}(H_{q})=\emptyset,
\end{eqnarray*}
see \cite[Theorem 3.1.4 in Chapter 1.3]{AG}.

\begin{rem} 
In \cite[Theorem 3.1.4 in Chapter 1.3]{AG}, 
it was shown that if $q\ge0$, $H_{q}$ has no eigenvalues 
and if $q<0$, $H_{q}$ has precisely one negative, 
simple eigenvalue. 
This fact is the difference between the repulsive and attractive cases. 
\end{rem}

Henceforth we consider the case $q\ge0$. 
To obtain the explicit formula for 
the unitary group $\{e^{-itH_{q}}\}_{t\in\rre}$, 
we introduce the Jost 
function associated to $H_{q}$. 

Let $f_{\pm}=f_{\pm}(x,\xi)$ be 
the unique solutions to the equation
\begin{eqnarray*}
H_{q}f=\frac12\xi^{2}f
\end{eqnarray*}
satisfying the asymptotic conditions 
\begin{eqnarray*}
f_{\pm}(x,\xi)-e^{\pm ix\xi}\to0,\quad as\ x\to\pm\infty.
\end{eqnarray*}
Then the elementary calculation yields 
\begin{eqnarray*}
f_{+}(x,\xi)&=&\left\{
\begin{array}{l}
\medskip\displaystyle{e^{ix\xi}\qquad\qquad\qquad\qquad\ \ \ \text{if}\ \ x\ge0,}\\
\displaystyle{\frac{1}{t_{q}(\xi)}e^{ix\xi}+
\frac{r_{q}(\xi)}{t_{q}(\xi)}e^{-ix\xi}\quad \text{if}\ \ x<0,}
\end{array}
\right.\\
f_{-}(x,\xi)&=&
\left\{
\begin{array}{l}
\medskip\displaystyle{\frac{1}{t_{q}(\xi)}e^{-ix\xi}+
\frac{r_{q}(\xi)}{t_{q}(\xi)}e^{ix\xi}\ \ \ \text{if}\ \ x\ge0,}\\
\displaystyle{e^{-ix\xi}\ \qquad\qquad\qquad\qquad\text{if}\ \ x<0,}
\end{array}
\right.
\end{eqnarray*}
where $t_{q}$ and $r_{q}$ are the transmission and reflection 
coefficients:
\begin{eqnarray}
t_{q}(\xi)=\frac{i\xi}{i\xi-q},\quad
r_{q}(\xi)=\frac{q}{i\xi-q}.\label{ed}
\end{eqnarray}
To obtain the the following representation formula 
for $e^{-itH_{q}}$, we introduce the following two 
operators:
\begin{eqnarray}
{{\mathcal L}}_{+}[\phi](x)
&=&
\phi(x)-q1_{-}(x)e^{qx}
\int_{x}^{-x}e^{q|y|}\phi(y)dy,
\label{K1}\\
{{\mathcal L}}_{-}[\phi](x)
&=&
\phi(x)-q1_{+}(x)e^{-qx}
\int_{-x}^{x}e^{q|y|}\phi(y)dy,
\label{K2}
\end{eqnarray}

\vskip3mm

\begin{prop}\label{l1}
Let $\phi\in L^{2}(\rre)$. Then we have
\begin{eqnarray}
(e^{-itH_{q}}\phi)(x)
&=&\left\{
\begin{array}{l}
\medskip\displaystyle{
\frac{e^{-i\pi/4}}{\sqrt{2\pi}}
t^{-1/2}
\int_{-\infty}^{\infty}e^{\frac{i(x-y)^{2}}{2t}}
{{\mathcal L}}_{+}[\phi](y)dy
\quad \text{if}\ \ x\ge0,}\\
\displaystyle{
\frac{e^{-i\pi/4}}{\sqrt{2\pi}}
t^{-1/2}
\int_{-\infty}^{\infty}e^{\frac{i(x-y)^{2}}{2t}}
{{\mathcal L}}_{-}[\phi](y)dy\quad \text{if}\ \ x<0,}
\end{array}
\right.\label{P}
\end{eqnarray}
\end{prop}

\vskip3mm

To prove this proposition, we show the following lemma.

\begin{lem}\label{l2}
Let $\phi\in L^{2}(\rre)$. Then we have 
\begin{eqnarray}
\lefteqn{(e^{-itH_{q}}\phi)(x)}\nonumber\\
&=&\left\{
\begin{array}{l}
\medskip
\displaystyle{
{{\mathcal F}}^{-1}[e^{-it\xi^{2}/2}
\{{{\mathcal F}}[\phi](\xi)
+r_{q}(\xi){{\mathcal F}}[1_{+}\phi](-\xi)
+r_{q}(\xi){{\mathcal F}}[1_{-}\phi](\xi)\}](x)
\quad \text{if}\ \ x\ge 0.}\\
\displaystyle{
{{\mathcal F}}^{-1}[e^{-it\xi^{2}/2}
\{{{\mathcal F}}[\phi](\xi)
+\overline{r_{q}(\xi)}{{\mathcal F}}[1_{-}\phi](-\xi)
+\overline{r_{q}(\xi)}{{\mathcal F}}[1_{+}\phi](\xi)\}](x)
\quad \text{if}\ \ x<0,}\\
\end{array}
\right.\nonumber\\
& &\label{Q}
\end{eqnarray}
\end{lem}

\vskip3mm
\noindent
{\bf Proof of Lemma \ref{l2}.} 
The explicit representation for the spectral projection
implies 
\begin{eqnarray}
\lefteqn{(e^{-itH_{q}}\phi)(x)}\nonumber\\
&=&\frac{1}{2\pi}
\int_{-\infty}^{\infty}
\left\{\int_{0}^{\infty}
e^{-\frac{i}{2}t\xi^{2}}
|t_{q}(\xi)|^{2}\right.
\nonumber\\
& &\qquad\qquad\left.\times\left(f_{+}(x,\xi)\overline{f_{+}(y,\xi)}
+f_{-}(x,\xi)\overline{f_{-}(y,\xi)}\right)d\xi
\right\}\phi(y)dy.\label{S}
\end{eqnarray}
We first consider the case $x\ge0$. 
In this case, we have 
\begin{eqnarray*}
\lefteqn{(e^{-itH_{q}}\phi)(x)}\\
&=&\frac{1}{2\pi}
\int_{0}^{\infty}
e^{-\frac{i}{2}t\xi^{2}}
|t_{q}(\xi)|^{2}f_{+}(x,\xi)
\left(\int_{0}^{\infty}
\overline{f_{+}(y,\xi)}
\phi(y)dy\right)d\xi\\
& &+\frac{1}{2\pi}
\int_{0}^{\infty}
e^{-\frac{i}{2}t\xi^{2}}
|t_{q}(\xi)|^{2}f_{+}(x,\xi)
\left(\int_{-\infty}^{0}
\overline{f_{+}(y,\xi)}
\phi(y)dy\right)d\xi\\
& &+\frac{1}{2\pi}
\int_{0}^{\infty}
e^{-\frac{i}{2}t\xi^{2}}
|t_{q}(\xi)|^{2}f_{-}(x,\xi)
\left(\int_{0}^{\infty}
\overline{f_{-}(y,\xi)}
\phi(y)dy\right)d\xi\\
& &+\frac{1}{2\pi}
\int_{0}^{\infty}
e^{-\frac{i}{2}t\xi^{2}}
|t_{q}(\xi)|^{2}f_{-}(x,\xi)
\left(\int_{-\infty}^{0}
\overline{f_{-}(y,\xi)}
\phi(y)dy\right)d\xi\\
&\equiv&F_{1}(t,x)+F_{2}(t,x)
+F_{3}(t,x)+F_{4}(t,x).
\end{eqnarray*}
From the definition of $f_{\pm}$, we have
\begin{eqnarray*}
F_{1}(t,x)
&=&\frac{1}{\sqrt{2\pi}}
\int_{0}^{\infty}
e^{ix\xi-\frac{i}{2}t\xi^{2}}
|t_{q}(\xi)|^{2}
{{\mathcal F}}[1_{+}\phi](\xi)d\xi,\\
F_{2}(t,x)
&=&\frac{1}{\sqrt{2\pi}}
\int_{0}^{\infty}
e^{ix\xi-\frac{i}{2}t\xi^{2}}
t_{q}(\xi)
{{\mathcal F}}[1_{-}\phi](\xi)d\xi\\
& &+\frac{1}{\sqrt{2\pi}}
\int_{0}^{\infty}
e^{ix\xi-\frac{i}{2}t\xi^{2}}
t_{q}(\xi)\overline{r_{q}(\xi)}
{{\mathcal F}}[1_{-}\phi](-\xi)d\xi\\
&\equiv&F_{21}(t,x)+F_{22}(t,x),\\
F_{3}(t,x)
&=&\frac{1}{\sqrt{2\pi}}
\int_{-\infty}^{0}
e^{ix\xi-\frac{i}{2}t\xi^{2}}
{{\mathcal F}}[1_{+}\phi](\xi)d\xi\\
& &+\frac{1}{\sqrt{2\pi}}
\int_{0}^{\infty}
e^{ix\xi-\frac{i}{2}t\xi^{2}}
r_{q}(\xi)
{{\mathcal F}}[1_{+}\phi](-\xi)d\xi\\
& &+\frac{1}{\sqrt{2\pi}}
\int_{-\infty}^{0}
e^{ix\xi-\frac{i}{2}t\xi^{2}}
r_{q}(\xi)
{{\mathcal F}}[1_{+}\phi](-\xi)d\xi\\
& &+\frac{1}{\sqrt{2\pi}}
\int_{0}^{\infty}
e^{ix\xi-\frac{i}{2}t\xi^{2}}
|r_{q}(\xi)|^{2}
{{\mathcal F}}[1_{+}\phi](\xi)d\xi\\
&\equiv&F_{31}(t,x)+F_{32}(t,x)+
F_{33}(t,x)+F_{34}(t,x),\\
F_{4}(t,x)
&=&\frac{1}{\sqrt{2\pi}}
\int_{-\infty}^{0}
e^{ix\xi-\frac{i}{2}t\xi^{2}}
t_{q}(\xi)
{{\mathcal F}}[1_{-}\phi](\xi)d\xi\\
& &
+\frac{1}{\sqrt{2\pi}}
\int_{0}^{\infty}
e^{ix\xi-\frac{i}{2}t\xi^{2}}
\overline{t_{q}(\xi)}r_{q}(\xi)
{{\mathcal F}}[1_{-}\phi](-\xi)d\xi\\
&\equiv&F_{41}(t,x)+F_{42}(t,x).
\end{eqnarray*}
Since $|t_{q}(\xi)|^{2}+|r_{q}(\xi)|^{2}=1$, we have 
\begin{eqnarray*}
F_{1}(t,x)+F_{31}(t,x)+F_{34}(t,x)=
{{\mathcal F}}^{-1}[e^{-it\xi^{2}/2}
{{\mathcal F}}[1_{+}\phi](\xi)](x).
\end{eqnarray*}
The identity $t_{q}(\xi)=1+r_{q}(\xi)$ yields 
\begin{eqnarray*}
\lefteqn{F_{21}(t,x)+F_{41}(t,x)}\\
&=&
{{\mathcal F}}^{-1}[e^{-it\xi^{2}/2}
{{\mathcal F}}[1_{-}\phi](\xi)](x)+
{{\mathcal F}}^{-1}[e^{-it\xi^{2}/2}
r_{q}(\xi){{\mathcal F}}[1_{-}\phi](\xi)](x).
\end{eqnarray*}
Using the relation $t_{q}(\xi)\overline{r_{q}(\xi)}+
\overline{t_{q}(\xi)}t_{q}(\xi)=0$, we obtain 
\begin{eqnarray*}
F_{22}(t,x)+F_{42}(t,x)=0.
\end{eqnarray*}
Combining the above identities with 
a trivial identity 
\begin{eqnarray*}
F_{32}(t,x)+F_{33}(t,x)=
{{\mathcal F}}^{-1}[e^{-it\xi^{2}/2}
r_{q}(\xi){{\mathcal F}}[1_{+}\phi](-\xi)](x),
\end{eqnarray*}
we have (\ref{Q}) with $x\ge0$. 

Next we consider the case $x<0$. 
Notice that $f_{+}(x,\xi)=f_{-}(-x,\xi)$. 
Then (\ref{S}) implies 
\begin{eqnarray*}
\lefteqn{(e^{-itH_{q}}\phi)(x)}\\
&=&\frac{1}{2\pi}
\int_{-\infty}^{\infty}
\biggl\{\int_{0}^{\infty}
e^{-\frac{i}{2}t\xi^{2}}
|t_{q}(\xi)|^{2}\\
& &\qquad\qquad
\times\left(f_{-}(-x,\xi)\overline{f_{-}(-y,\xi)}
+f_{+}(-x,\xi)\overline{f_{+}(-y,\xi)}\right)d\xi
\biggl\}\phi(y)dy\\
&=&\frac{1}{2\pi}
\int_{-\infty}^{\infty}
\biggl\{\int_{0}^{\infty}
e^{-\frac{i}{2}t\xi^{2}}
|t_{q}(\xi)|^{2}.\\
& &\qquad\qquad
\times\left(f_{-}(-x,\xi)\overline{f_{-}(y,\xi)}
+f_{+}(-x,\xi)\overline{f_{+}(y,\xi)}\right)d\xi
\biggl\}\phi(-y)dy\\
&=&(e^{-itH_{q}}\phi(-\cdot))(-x).
\end{eqnarray*}
Hence the identity 
(\ref{Q}) with $x<0$ follows from (\ref{Q}) with 
$x\ge0$. This completes the proof of Lemma \ref{l2}.
$\qquad\qed$

\vskip3mm
\noindent
{\bf Proof of Proposition \ref{l1}.}
We consider the case $x\ge0$ only because the identity 
(\ref{P}) with the case $x<0$ 
follows from (\ref{P}) with the case $x\ge0$ through the relation 
$(e^{-itH_{q}}\phi)(x)=(e^{-itH_{q}}\phi(-\cdot))(-x)$. 
By Lemma \ref{l2}, we obtain
\begin{eqnarray*}
\lefteqn{(e^{-itH_{q}}\phi)(x)}\\
&=&(\frac{1}{\sqrt{2\pi}}
{{\mathcal F}}^{-1}[e^{-it\xi^{2}/2}]\ast \phi)(x)
+\frac{1}{2\pi}({{\mathcal F}}^{-1}[e^{-it\xi^{2}/2}]
\ast{{\mathcal F}}^{-1}[r_{q}]\ast(1_{+}\phi)(-\cdot))(x)\\
& &+\frac{1}{2\pi}({{\mathcal F}}^{-1}[e^{-it\xi^{2}/2}]
\ast{{\mathcal F}}^{-1}[r_{q}]\ast(1_{-}\phi))(x).
\end{eqnarray*}
Since ${{\mathcal F}}^{-1}[r_{q}](x)=-\sqrt{2\pi}q1_{-}(x)e^{qx}$, 
we obtain
\begin{eqnarray*}
({{\mathcal F}}^{-1}[r_{q}]\ast(1_{+}\phi)(-\cdot))(x)
&=&-\sqrt{2\pi}q1_{-}(x)e^{qx}
\int_{0}^{-x}e^{qy}\phi(y)dy,\\
({{\mathcal F}}^{-1}[r_{q}]\ast1_{-}\phi)(x)
&=&-\sqrt{2\pi}q1_{-}(x)e^{qx}
\int^{0}_{x}e^{-qy}\phi(y)dy.
\end{eqnarray*}
Combining the above identities with the well-known 
identity 
\begin{eqnarray*}
{{\mathcal F}}^{-1}[e^{-it\xi^{2}/2}](x)
=e^{-i\pi/4}t^{-1/2}e^{\frac{ix^{2}}{2t}}, 
\end{eqnarray*}
we obtain (\ref{P}). 
$\qquad\qed$

\vskip3mm

Next lemma is concerned with the Strichartz 
estimates for the group $e^{-itH_{q}}$ 
which plays an important role to construct 
the solution to (\ref{NLS}) in the finite 
interval $[0,T)$ (see Lemma \ref{wp} below) 
and the infinite interval $[T,\infty)$ 
(see Section 4). 

\begin{lem}\label{S} Let 
$2/q_j+1/r_j=1/2$, $4\le q_j\le\infty$ and $j=1,2$. 
For any (possibly unbounded) interval $I$ and 
for any $s\in\overline{I}$, we have
\begin{eqnarray}
\|\int_{s}^{t}e^{-i(t-\tau)H_{q}}
F(\tau)d\tau\|_{L^{q_1}(I;L_x^{r_1})}
\lesssim\|F\|_{L^{q_2'}(I;L_x^{r_2'})},
\label{U}
\end{eqnarray}
where $p'$ is the H\"{o}lder conjugate exponent of $p$, 
and the implicit constant depends only on $p$.
\end{lem}

\vskip3mm
\noindent
{\bf Proof of Lemma \ref{S}.} Although Lemma \ref{S} is proved in 
\cite[Proposition 2.2]{HMZ2}, we give the proof of this lemma for the reader's 
convenience. From the $L^{2}$ unitarity of $e^{-itH_{q}}$, 
we have 
\begin{eqnarray*}
\|e^{-itH_{q}}\phi\|_{L_{x}^{2}(\rre)}&=&
\|\phi\|_{L_{x}^{2}(\rre)}.
\end{eqnarray*}
From Proposition \ref{l1} (\ref{P}) and the fact that 
${{\mathcal L}}_{\pm}:L^{1}\to L^{1}$ is bounded, we have
\begin{eqnarray*}
\|e^{-itH_{q}}\phi\|_{L_{x}^{\infty}(\rre)}&\lesssim&
t^{-1/2}(\|{{\mathcal L}}_{+}[\phi]\|_{L^{1}_{x}}
+\|{{\mathcal L}}_{-}[\phi]\|_{L^{1}_{x}})\\
&\lesssim&t^{-1/2}\|\phi\|_{L^{1}_{x}}. 
\end{eqnarray*}
Combining the above two boundness for $e^{-itH_{q}}$ 
with $TT^{\ast}$ argument 
(see Yajima \cite{Y} and Keel-Tao\cite{KT} 
for instance) we obtain (\ref{U}).
$\qquad\qed$

\vskip3mm

We next state the lemma concerning the global well-posedness 
for (\ref{NLS}) in $L^{2}(\rre)$. 

\begin{lem}\label{wp} 
For any $\phi\in L^{2}(\rre)$, there exists a unique solution 
$u\in C(\rre;L^{2}(\rre))\cap L^{q}_{loc}(\rre;L_{x}^{r}(\rre))$ 
to (\ref{NLS}) with $u(0,x)=\phi(x)$, 
where $(q,r)$ satisfies 
$2/q+1/r=1/2$ and $4\le q\le\infty$. 
\end{lem}

\vskip3mm
\noindent
{\bf Proof of Lemma \ref{wp}.} The proof follows from 
the the Strichartz estimate (Lemma \ref{S} (\ref{U}) 
with $I=[0,T)$ and $s=0$) and 
the fixed point argument, see Tsutsumi \cite{T} 
for the detail. $\qquad\qed$

\section{Distorted Fourier Transform} \label{sec:dist}

In this section, we introduce a distorted Fourier transform 
associated to $H_{q}$. 

Let us review that the usual Fourier transform and its inverse transform 
are defined by
\begin{eqnarray*}
{{\mathcal F}}[\phi](\xi)
&=&\frac{1}{\sqrt{2\pi}}
\int_{-\infty}^{+\infty}e^{-iy\xi}\phi(y)dy,\\
{{\mathcal F}}^{-1}[\phi](x)
&=&\frac{1}{\sqrt{2\pi}}
\int_{-\infty}^{+\infty}e^{ix\xi}\phi(\xi)d\xi.
\end{eqnarray*}
Let $t_{q}$ and $r_{q}$ are defined by 
(\ref{ed}). 
We introduce 
\begin{eqnarray*}
\Psi(x,\xi)=\left\{
\begin{array}{l}
\medskip\displaystyle{t_{q}(\xi)f_{+}(x,\xi)\qquad\ 
\ \ \ \ \text{if}\ \xi\ge0,}\\
\displaystyle{t_{q}(-\xi)f_{-}(x,-\xi)\qquad 
\text{if}\  \xi<0.}
\end{array}
\right.
\end{eqnarray*}
Then the distorted Fourier transform and its inverse transform 
for $H_{q}$ are defined by
\begin{eqnarray*}
{{\mathcal F}}_{q}[\phi](\xi)
&=&\frac{1}{\sqrt{2\pi}}
\int_{-\infty}^{+\infty}\Psi(y,-\xi)\phi(y)dy,\\
{{\mathcal F}}_{q}^{-1}[\phi](x)
&=&\frac{1}{\sqrt{2\pi}}
\int_{-\infty}^{+\infty}\overline{\Psi(x,-\xi)}\phi(\xi)d\xi.
\end{eqnarray*}
By the definition, the relations between the usual and 
distorted Fourier transforms are as follows:
\begin{eqnarray}
{{\mathcal F}}_{q}[\phi](\xi)
&=&\left\{
\begin{array}{l}
\medskip\displaystyle{
{{\mathcal F}}[\phi](\xi)
+r_{q}(\xi){{\mathcal F}}[1_{+}\phi](-\xi)
+r_{q}(\xi){{\mathcal F}}[1_{-}\phi](\xi)
\quad\text{if}\ \xi\ge0,}\\
\displaystyle{
{{\mathcal F}}[\phi](\xi)
+\overline{r_{q}(\xi)}{{\mathcal F}}[1_{-}\phi](-\xi)
+\overline{r_{q}(\xi)}{{\mathcal F}}[1_{+}\phi](\xi)
\quad\text{if}\ \xi<0,}
\end{array}
\right.\label{F}
\end{eqnarray}
\begin{eqnarray}
{{\mathcal F}}_{q}^{-1}[\phi](x)
&=&\left\{
\begin{array}{l}
\medskip\displaystyle{
{{\mathcal F}}^{-1}[\phi](x)
+{{\mathcal F}}^{-1}[1_{+}\overline{r_{q}}\phi](-x)
+{{\mathcal F}}^{-1}[1_{-}r_{q}\phi](x)
\quad\text{if}\ x\ge0,}\\
\displaystyle{
{{\mathcal F}}^{-1}[\phi](x)
+{{\mathcal F}}^{-1}[1_{-}r_{q}\phi](-x)
+{{\mathcal F}}^{-1}[1_{+}\overline{r_{q}}\phi](x)
\quad\text{if}\ x<0,}
\end{array}
\right.\label{G}
\end{eqnarray}
Furthermore, we have the another representation for 
the distorted Fourier transforms:
\begin{eqnarray}
{{\mathcal F}}_{q}[\phi](\xi)
&=&\left\{
\begin{array}{l}
\medskip\displaystyle{
{{\mathcal F}}[\phi(x)-q1_{-}(x)e^{qx}
\int_{x}^{-x}e^{q|y|}\phi(y)dy](\xi)
\ \quad\text{if}\ \xi\ge0,}
\nonumber\\
\displaystyle{
{{\mathcal F}}[\phi(x)-q1_{+}(x)e^{-qx}
\int_{-x}^{x}e^{q|y|}\phi(y)dy](\xi)
\quad\text{if}\ \xi<0,}
\end{array}
\right.\\
&=&\left\{
\begin{array}{l}
\medskip\displaystyle{
{{\mathcal F}}[{{\mathcal L}}_{+}[\phi](x)](\xi)
\ \quad\text{if}\ \xi\ge0,}\\
\displaystyle{
{{\mathcal F}}[{{\mathcal L}}_{-}[\phi](x)](\xi)
\quad\text{if}\ \xi<0,}
\end{array}
\right.\label{I}
\end{eqnarray}
\begin{eqnarray}
{{\mathcal F}}_{q}^{-1}[\phi](x)
&=&\left\{
\begin{array}{l}
\medskip\displaystyle{
{{\mathcal F}}^{-1}[\phi](x)
-qe^{qx}\int_{-\infty}^{-x}e^{qy}
{{\mathcal F}}^{-1}[1_{+}\phi](y)dy}\\
\medskip\displaystyle{
\qquad\qquad-qe^{qx}\int_{x}^{+\infty}e^{-qy}
{{\mathcal F}}^{-1}[1_{-}\phi](y)dy
\quad\text{if}\ x\ge0,}\\
\medskip\displaystyle{
{{\mathcal F}}^{-1}[\phi](x)
-qe^{-qx}\int_{-\infty}^{x}e^{qy}{{\mathcal F}}^{-1}[1_{+}\phi](y)dy}\\
\medskip\displaystyle{
\qquad\qquad-qe^{-qx}\int^{+\infty}_{-x}e^{-qy}
{{\mathcal F}}^{-1}[1_{-}\phi](y)dy
\quad\text{if}\ x<0.}
\end{array}
\right.\label{J}
\end{eqnarray}
From (\ref{I}) and (\ref{J}), we see that 
\begin{eqnarray*}
{{\mathcal F}}_{q}{{\mathcal F}}_{q}^{-1}
={{\mathcal F}}_{q}^{-1}{{\mathcal F}}_{q}=I\qquad 
\text{in}\ L^{2}(\rre)
\end{eqnarray*}
and
\begin{eqnarray*}
\int_{\rre}
{{\mathcal F}}_{q}[\phi](x)
\overline{\varphi(x)}dx
=\int_{\rre}\phi(x)
\overline{{{\mathcal F}}_{q}^{-1}[\varphi](x)}dx
\end{eqnarray*}
for any $\phi,\varphi\in L^{2}(\rre)$. 
From the above identities, 
we have the $L^{2}$ isometry of ${{\mathcal F}}_{q}$ and 
${{\mathcal F}}_{q}^{-1}$:
\begin{eqnarray*}
\|{{\mathcal F}}_{q}[\phi]\|_{L^{2}_{\xi}}
=\|\phi\|_{L^{2}_{x}},\qquad
\|{{\mathcal F}}_{q}^{-1}[\varphi]\|_{L^{2}_{x}}
=\|\varphi\|_{L^{2}_{\xi}}.
\end{eqnarray*}

\begin{rem} We can express 
the unitary group $e^{-itH_{q}}$ 
in terms of the distorted Fourier transform. 
Indeed, we obtain
\begin{eqnarray*}
(e^{-itH_{q}}\phi)(x)
=\overline{{{\mathcal F}}_{q}^{-1}
\left[e^{\frac{i}{2}t\xi^{2}}
{{\mathcal F}}_{q}\left[\overline{\phi}
\right](\xi)\right](x)}
\end{eqnarray*}
for any $\phi\in L^{2}(\rre)$.
\end{rem}
\vskip3mm

Next proposition is concerned with the asymptotic behavior of 
the unitary group $e^{-itH_{q}}$.

\begin{prop}\label{l3} 
Let $\phi\in{{\mathcal S}}'(\rre)$ 
satisfy $\phi$, $\pt_{\xi}
\phi\in L_{\xi}^{2}$ and $\phi(0)=0$. Let $2\le p\le\infty$ 
and $\beta$ satisfy 
$0\le\beta\le1$ for $p=2$ and 
$0\le\beta<(p+2)/(2p)$ for $2<p\le\infty$. 
Then we have
\begin{eqnarray*}
(e^{-itH_{q}}{{\mathcal F}}_{q}^{-1}\phi)(x)
=t^{-1/2}\phi(\frac{x}{t})
e^{\frac{ix^{2}}{2t}-\frac{i\pi}{4}}+R(t,x),
\end{eqnarray*}
where $R$ satisfies
\begin{eqnarray}
\|R(t)\|_{L^{p}_{x}}
\lesssim 
t^{-1/2+1/p-\beta/2}
(\|\phi\|_{L_{\xi}^{2}}
+\|\pt_{\xi}\phi\|_{L_{\xi}^{2}})
\label{r}
\end{eqnarray}
for $t>0$.
\end{prop}

\vskip3mm

To prove Proposition \ref{l3}, we show the following 
identities.

\begin{lem} \label{ly} 
Let ${{\mathcal L}}_{\pm}$ be defined by 
(\ref{K1}) and (\ref{K2}). Then we have
\begin{eqnarray}
{{\mathcal L}}_{+}
[{{\mathcal F}}_{q}^{-1}[\phi]](y)&=&
{{\mathcal F}}_{q,+}^{-1}[\phi](y),
\label{e1}\\
{{\mathcal L}}_{-}
[{{\mathcal F}}_{q}^{-1}[\phi]](y)&=&
{{\mathcal F}}_{q,-}^{-1}[\phi](y)
\label{e2}
\end{eqnarray}
for any $y\in\rre$, 
where ${{\mathcal F}}_{q,\pm}^{-1}$ are defined by
\begin{eqnarray*}
{{\mathcal F}}_{q,+}^{-1}[\phi](y)
&=&{{\mathcal F}}^{-1}[\phi](x)
-qe^{qx}\int_{-\infty}^{-x}e^{qy}
{{\mathcal F}}^{-1}[1_{+}\phi](y)dy\\
& &\qquad-qe^{qx}\int_{x}^{+\infty}e^{-qy}
{{\mathcal F}}^{-1}[1_{-}\phi](y)dy,\\
{{\mathcal F}}_{q,-}^{-1}[\phi](y)
&=&
{{\mathcal F}}^{-1}[\phi](x)
-qe^{-qx}\int_{-\infty}^{x}e^{qy}
{{\mathcal F}}^{-1}[1_{+}\phi](y)dy\\
& &\qquad
-qe^{-qx}\int^{+\infty}_{-x}e^{-qy}
{{\mathcal F}}^{-1}[1_{-}\phi](y)dy.
\end{eqnarray*}
\end{lem}

\vskip3mm

\begin{rem}
The operators ${{\mathcal F}}_{q,\pm}^{-1}$ are 
appeared in the representation formula (\ref{J}) 
for ${{\mathcal F}}_{q}^{-1}$. 
\end{rem}

\vskip3mm

We shall prove Lemma \ref{ly} in Appendix below. 

\vskip3mm
\noindent
{\bf Proof of Proposition \ref{l3}.} 
We consider the case $x\ge0$ only. 
From Proposition \ref{l1} (\ref{P}) we have
\begin{eqnarray*}
(e^{-itH_{q}}{{\mathcal F}}_{q}^{-1}\phi)(x)
&=&
\frac{e^{-i\pi/4}}{\sqrt{2\pi}}
t^{-1/2}e^{\frac{ix^{2}}{2t}}
\int_{-\infty}^{+\infty}e^{-i\frac{x}{t}y}
e^{\frac{iy^{2}}{2t}}
{{\mathcal L}}_{+}[{{\mathcal F}}_{q}^{-1}[\phi]](y)dy\\
&=&
\frac{e^{-i\pi/4}}{\sqrt{2\pi}}
t^{-1/2}e^{\frac{ix^{2}}{2t}}
\int_{-\infty}^{+\infty}e^{-i\frac{x}{t}y}
{{\mathcal L}}_{+}[{{\mathcal F}}_{q}^{-1}[\phi]](y)dy\\
& &+
\frac{e^{-i\pi/4}}{\sqrt{2\pi}}
t^{-1/2}e^{\frac{ix^{2}}{2t}}
\int_{-\infty}^{+\infty}e^{-i\frac{x}{t}y}
(e^{\frac{iy^{2}}{2t}}-1)
{{\mathcal L}}_{+}[{{\mathcal F}}_{q}^{-1}[\phi]](y)dy\\
&\equiv&L(t,x)+R(t,x).
\end{eqnarray*}
From the representation formula (\ref{I}) for 
${{\mathcal F}}_{q}$, we see 
\begin{eqnarray*}
L(t,x)
=t^{-1/2}\phi(\frac{x}{t})
e^{\frac{ix^{2}}{2t}-\frac{i\pi}{4}}.
\end{eqnarray*}
For $R$, we evaluate as 
\begin{eqnarray*}
\|R(t)\|_{L_{x}^{p}(0,\infty)}
&\lesssim&t^{-1/2+1/p}
\|(e^{\frac{iy^{2}}{2t}}-1)
{{\mathcal L}}_{+}[{{\mathcal F}}_{q}^{-1}[\phi]](y)
\|_{L_{y}^{p'}(\rre)}\\
&\lesssim&t^{-1/2+1/p-\beta/2}
\||y|^{\beta}
{{\mathcal L}}_{+}[{{\mathcal F}}_{q}^{-1}[\phi]](y)
\|_{L_{y}^{p'}(\rre)}\\
&\lesssim&t^{-1/2+1/p-\beta/2}
\|\langle y\rangle^{\beta+\gamma}
{{\mathcal L}}_{+}[{{\mathcal F}}_{q}^{-1}[\phi]](y)
\|_{L_{y}^{2}(\rre)},
\end{eqnarray*}
where $2\le p\le\infty$, 
$0\le\beta\le1$ and $\gamma$ satisfies 
$\gamma=0$ for $p=2$ and 
$\gamma>(p-2)/(2p)$ for $2<p\le\infty$. 
We choose $\beta$ so that $\beta+\gamma=1$. 
Using Lemma \ref{ly} and the representation (\ref{G}) 
for ${{\mathcal F}}_{q,+}^{-1}$, 
we obtain
\begin{eqnarray*}
\lefteqn{\|\langle y\rangle
{{\mathcal L}}_{+}[\phi](y)
\|_{L_{y}^{2}(\rre)}}\\
&\le&
\|{{\mathcal F}}^{-1}[\phi](y)
\|_{L_{y}^{2}(\rre)}
+
\|y{{\mathcal F}}^{-1}[\phi](y)
\|_{L_{y}^{2}(\rre)}\\
& &
+\|{{\mathcal F}}^{-1}[1_{+}\overline{r_{q}}
\phi](-y)\|_{L_{y}^{2}(\rre)}
+\|y{{\mathcal F}}^{-1}[1_{+}\overline{r_{q}}
\phi](-y)\|_{L_{y}^{2}(\rre)}
\\
& &
+\|{{\mathcal F}}^{-1}[1_{-}r_{q}\phi](y)
\|_{L_{y}^{2}(\rre)}
+\|y{{\mathcal F}}^{-1}[1_{-}r_{q}\phi](y)
\|_{L_{y}^{2}(\rre)}\\
&\equiv&I_{1}+I_{2}+I_{3}+I_{4}+I_{5}+I_{6}.
\end{eqnarray*}
Since $|r_{q}(\xi)|+|r_{q}'(\xi)|\lesssim1$, we see that 
\begin{eqnarray*}
I_{1}+I_{3}+I_{5}\lesssim
\|{{\mathcal F}}_{q}[\phi]\|_{L_{\xi}^{2}}.
\end{eqnarray*}
Notice that if $\varphi(0)=0$, then 
$y{{\mathcal F}}^{-1}[1_{\pm}\varphi](y)
=i{{\mathcal F}}^{-1}[1_{\pm}\pt_{\xi}\varphi](y)$. 
Hence we have
\begin{eqnarray*}
I_{2}&=&\|{{\mathcal F}}^{-1}[\pt_{\xi}\phi]
\|_{L_{y}^{2}(\rre)}
\lesssim\|\pt_{\xi}\phi\|_{L_{\xi}^{2}},\\
I_{4}&=&\|{{\mathcal F}}^{-1}[1_{+}\pt_{\xi}
(\overline{r}_{q}\phi)]
\|_{L_{y}^{2}(\rre)}
\lesssim\|\phi\|_{L_{\xi}^{2}}
+\|\pt_{\xi}\phi\|_{L_{\xi}^{2}},\\
I_{6}&=&\|{{\mathcal F}}^{-1}[1_{-}\pt_{\xi}
(r_{q}\phi)]
\|_{L_{y}^{2}(\rre)}
\lesssim\|\phi\|_{L_{\xi}^{2}}
+\|\pt_{\xi}\phi\|_{L_{\xi}^{2}}.
\end{eqnarray*}
Therefore we obtain (\ref{r}) for the case $x\ge0$. 
By a similar way as above, 
we obtain (\ref{r}) for the case $x<0$ 
which completes the proof of Proposition \ref{l3}. 
$\qquad\qed$

\section{Construction of Modified Wave Operators} \label{sec:wave}

In this section we prove Theorem \ref{main}. 
We first 
rewrite (\ref{NLS}) as the integral equation. Using 
Proposition \ref{l3}, we have
\begin{eqnarray}
u(t)-u_{ap}(t)=u(t)-e^{-itH_{q}}{{\mathcal F}}_{q}^{-1}w
+R_1(t),\label{71}
\end{eqnarray}
where 
\begin{eqnarray}
w(t,\xi)={{\mathcal F}}_{q}[\phi_{+}](\xi)
e^{-i\lambda|{{\mathcal F}}_{q}[\phi_{+}](\xi)|^{2}\log t}, \label{72}
\end{eqnarray}
and $R_1$ satisfies 
\begin{eqnarray*}
\lefteqn{\|R_1(t)\|_{L^{\infty}(t,\infty;L_x^2)}}\\
&\lesssim&
t^{-\beta/2}(\|{{\mathcal F}}_{q}[\phi_+]\|_{L_{\xi}^{2}}
+\|\pt_{\xi}{{\mathcal F}}_{q}[\phi_+]\|_{L_{\xi}^{2}})
(1+\|{{\mathcal F}}_{q}[\phi_+]\|_{L_{\xi}^{2}}
+\|\pt_{\xi}{{\mathcal F}}_{q}[\phi_+]\|_{L_{\xi}^{2}})^{2}
\end{eqnarray*}
for $0\le\beta\le1$, where we used the fact 
that ${{\mathcal F}}_{q}[\phi_{+}](0)=0$. 
From the representation formula (\ref{F}) for ${{\mathcal F}}_{q}$, 
we see that
\begin{eqnarray*}
\|{{\mathcal F}}_{q}[\phi_{+}]\|_{L_{\xi}^{2}}
=\|\phi_{+}\|_{L_{x}^{2}},\qquad\|\pt_{\xi}
{{\mathcal F}}_{q}[\phi_{+}]\|_{L_{\xi}^{2}}
\lesssim\|\phi_{+}\|_{H_{x}^{0,1}}. 
\end{eqnarray*}
Hence we have
\begin{eqnarray}
\|R_1(t)\|_{L^{\infty}(t,\infty;L_x^2)}
\lesssim
t^{-\beta/2}\|\phi_+\|_{H_{x}^{0,1}}(1+\|\phi_+\|_{H_{x}^{0,1}})^{2}.
\label{77}
\end{eqnarray}
By a similar way we obtain 
\begin{eqnarray}
\|R_1(t)\|_{L^4(t,\infty;L_x^{\infty})}
&\lesssim&
t^{-\gamma/2-1/4}\|\phi_+\|_{H_{x}^{0,1}}
(1+\|\phi_+\|_{H_{x}^{0,1}})^{2},
\label{78}
\end{eqnarray}
where $0\le\gamma<1/2$. Eq. (\ref{71}) is rewritten as follows:
\begin{eqnarray}
u(t)-u_{ap}(t)=e^{-itH_{q}}{{\mathcal F
}}_{q}^{-1}[{{\mathcal F}}_{q}e^{itH_{q}}u-w]
+R_1(t).\label{73}
\end{eqnarray}
From (\ref{NLS}) and (\ref{72}), we obtain
\begin{eqnarray}
i\pt_t({{\mathcal F}}_{q}e^{itH_{q}}u)&=&
\lambda{{\mathcal F}}_{q}e^{itH_{q}}|u|^2u,\label{74}\\
i\pt_tw&=&\lambda t^{-1}
|w|^{2}w.\label{75}
\end{eqnarray}
Subtracting (\ref{75}) from (\ref{74}), we have
\begin{eqnarray}
i\pt_t({{{{\mathcal F}}}}_{q}e^{itH_{q}}u-w)
=\lambda{{{{\mathcal F}}}}_{q}e^{itH_{q}}|u|^2u
-\lambda t^{-1}|w|^{2}w.
\label{92}
\end{eqnarray}
Proposition \ref{l3} yields 
\begin{eqnarray}
e^{-itH_{q}}{{\mathcal F}}_{q}^{-1}
[\lambda t^{-1}|w|^{2}w]
=\lambda|u_{ap}|^2u_{ap}+R_2(t,x),
\label{87}
\end{eqnarray}
where $R_2$ satisfies 
\begin{eqnarray}
\|R_2\|_{L^1(t,\infty;L_x^2)}\lesssim
t^{-\beta/2}\|\phi_+\|_{H_{x}^{0,1}}(1+\|\phi_+\|_{H^{0,1}})^4,\label{79}
\end{eqnarray}
where $0\le\beta\le1$. 
Substituting (\ref{87}) into (\ref{92}), we obtain 
\begin{eqnarray*}
i\pt_t({{{{\mathcal F}}}}_{q}e^{itH_{q}}u-w)
=\lambda{{{{\mathcal F}}}}_{q}e^{itH_{q}}[|u|^2u-|u_{ap}|^2u_{ap}]
-{{{{\mathcal F}}}}_{q}e^{itH_{q}}R_2.
\end{eqnarray*}
Integrating the above equation in $t$,
we have
\begin{eqnarray}
u(t)-
e^{-itH_{q}}{{\mathcal F}}_{q}^{-1}
w&=&i\lambda\int_t^{+\infty}
e^{-i(t-\tau)H_{q}}(|u|^2u-|u_{ap}|^2u_{ap})(\tau)d\tau\nonumber\\
& &-i\int_t^{+\infty}e^{-i(t-\tau)H_{q}}R_2(\tau)d\tau.\label{80}
\end{eqnarray}
Combining (\ref{71}) with (\ref{80}), we obtain the following 
integral equation:
\begin{eqnarray}
u(t)-u_{ap}(t)&=&
i\lambda\int_t^{+\infty}
e^{-i(t-\tau)H_{q}}[|u|^2u-|u_{ap}|^2u_{ap}](\tau)d\tau\nonumber\\
& &+R_1(t)-i\int_t^{+\infty}e^{-i(t-\tau)
H_{q}}R_2(\tau)d\tau.\label{81}
\end{eqnarray}

To show the existence of $u$ satisfying (\ref{81}), we 
shall prove that the map
\begin{eqnarray*}
\Phi[u](t)&=&u_{ap}(t)+
i\lambda\int_t^{+\infty}
e^{-i(t-\tau)H_{q}}[|u|^2u-|u_{ap}|^2u_{ap}](\tau)d\tau\\
& &+R_1(t)-i\int_t^{+\infty}e^{-i(t-\tau)H_{q}}R_2(\tau)d\tau
\end{eqnarray*}
is a contraction on
\begin{eqnarray*}
{{{\bf X}}}_{r,T}&=&
\{u\in C([T,\infty);L^2(\rre))\cap L_{loc}^4(T,\infty;L^{\infty}(\rre));
\|u-u_{ap}\|_{{{\bf X}}_T}\le r\},\\
\|v\|_{{{\bf X}}_T}
&=&\sup_{t\ge T}
t^{\alpha}(\|v\|_{L^{\infty}(t,\infty;L_x^2)}
+\|v\|_{L^4(t,\infty;L_x^{\infty})})
\end{eqnarray*}
for some $r>0$ and $T>0$. 

Let $v(t)=u(t)-u_{ap}(t)$ and $\|\phi_{+}\|_{H^{0,1}}\le r$. Then 
for any $v\in {{\bf X}}_T$ satisfying $\|v\|_{{{\bf X}}_T}\le r$, 
we obtain
\begin{eqnarray*}
\lefteqn{\Phi[u](t)-u_{ap}(t)}\\
&=&
i\lambda\int_t^{+\infty}
e^{-i(t-\tau)H_{q}}[|v|^2v+2u_{ap}|v|^2+\overline{u}_{ap}
v^2+2|u_{ap}|^2v+u_{ap}^2\overline{v}](\tau)d\tau\\
& &+R_1(t,x)-i\int_t^{+\infty}e^{-i(t-\tau)H_{q}}R_2(\tau)d\tau.
\end{eqnarray*}
Strichartz estimate (Lemma \ref{S} (\ref{U}) 
with $I=[t,\infty)$ and $s=\infty$) implies
\begin{eqnarray}
\lefteqn{\|\Phi[u]-u_{ap}\|_{L^{\infty}(t,\infty;L_x^2)}
+\|\Phi[u]-u_{ap}\|_{L^4(t,\infty;L_x^{\infty})}}\nonumber\\
&\lesssim&(\|F_1\|_{L^{4/3}(t,\infty;L_x^1)}
+\|F_2\|_{L^1(t,\infty;L_x^2)}
+\|F_3\|_{L^1(t,\infty;L_x^2))}\nonumber\\
& &+(\|R_1\|_{L^{\infty}(t,\infty;L_x^2)}
+\|R_1\|_{L^4(t,\infty;L_x^{\infty})}+
\|R_2\|_{L^1(t,\infty;L_x^2)}),\label{82}
\end{eqnarray}
where $F_1=|v|^2v$, $F_2=2u_{ap}|v|^2+\overline{u}_{ap}v^2$ and 
$F_3=2|u_{ap}|^2v+u_{ap}^2\overline{v}$. By the H\"{o}lder inequality, 
\begin{eqnarray*}
\|F_1\|_{L^{4/3}(t,\infty;L_x^1)}
&\le&\|\|v\|_{L_x^2}^2\|v\|_{L_x^{\infty}}\|_{L^{4/3}(t,\infty)}
\nonumber\\
&\lesssim&
r^2\|t^{-2\alpha}\|v\|_{L_x^{\infty}}\|_{L^{4/3}(t,\infty)}\nonumber\\
&\lesssim&r^2\|t^{-2\alpha}\|_{L^2(t,\infty)}\|v\|_{L^4(t,\infty;L_x^{\infty})}
\nonumber\\
&\lesssim&t^{-3\alpha+1/2}r^3,\nonumber\\
\|F_2\|_{L^1(t,\infty;L_x^2)}
&\le&\|\|u_{ap}\|_{L_x^{\infty}}
\|v\|_{L_x^2}\|v\|_{L_x^{\infty}}\|_{L^1(t,\infty)}
\nonumber\\
&\lesssim& 
r^2\|t^{-\alpha-1/2}\|v\|_{L_x^{\infty}}\|_{L^1(t,\infty)}\\
&\lesssim&r^2\|t^{-\alpha-1/2}\|_{L^{4/3}(t,\infty)}
\|v\|_{L^4(t,\infty;L_x^{\infty})}
\nonumber\\
&\lesssim& t^{-2\alpha+1/4}r^3,\nonumber\\
\|F_3\|_{L^1(t,\infty;L_x^2)}
&\le&\|\|u_{ap}\|_{L_x^{\infty}}^2
\|v\|_{L_x^2}\|_{L^1(t,\infty)}
\nonumber\\
&\lesssim& r^3\|t^{-\alpha-1}\|_{L^1(t,\infty)}
\nonumber\\
&\lesssim& t^{-\alpha}r^3.
\end{eqnarray*}
Substituting above three inequalities, 
(\ref{77}), (\ref{78}) and (\ref{79}) 
into (\ref{82}), we have
\begin{eqnarray*}
\lefteqn{\|\Phi[u]-u_{ap}\|_{{{\bf X}}_T}}\\
&\le&C(T^{-2\alpha+1/2}+1)r^3+C
(T^{\alpha-\beta/2}+T^{\alpha-\gamma/2-1/4})r(1+r^4).
\end{eqnarray*}
Choosing $\alpha,\beta$ and $\gamma$ 
so that $1/4<\alpha<1/2$, $2\alpha<\beta<1$ 
and $2(\alpha-1/4)<\gamma<1/2$, 
and taking $T$ large enough and $r$ sufficiently small, 
we guarantee that $\Phi$ is the map onto ${{\bf X}}_{r,T}$. 
By a similar way we can conclude that  $\Phi$ is the contraction 
map on ${{\bf X}}_{r,T}$. Therefore Banach fixed point theorem 
yields that $\Phi$ has a unique fixed point in ${{\bf X}}_{r,T}$ 
which is the solution to the final state problem (\ref{NLS})-(\ref{ac}). 

Next, we show that the solution to (\ref{NLS}) with   
finite ${{\bf X}}_T$ norm is unique. Let 
$u$ and $v$ be two solutions satisfying 
$\|u\|_{{{\bf X}}_T}<\infty$ and $\|v\|_{{{\bf X}}_T}<\infty$. 
We put $t_1=\inf\{t\in[T,\infty);u(s)=v(s)$ for any $s\in[t,\infty)\}$ 
and $R=\max\{\|u\|_{{{\bf X}}_T},\|v\|_{{{\bf X}}_T}\}$. 
If $t_1=T$, then $u(t)=v(t)$ on $[T,\infty)$ which is desired 
result. If $T<t_1$, as in (\ref{82}) 
by the Strichartz inequality (Lemma \ref{S}), 
we have
\begin{eqnarray*}
\|u-v\|_{L^4(t_0,t_1,L_x^{\infty})}
\le CR^2(t_0^{1-4\alpha}-t_1^{1-4\alpha})^{1/2}
\|u-v\|_{L^4(t_0,t_1,L_x^{\infty})},
\end{eqnarray*}
for $t_0\in[T,t_1)$. Since $-2\alpha+1/2<0$, 
we can choose $t_0\in[T,t_1)$ so that 
$CR^2(t_0^{1-4\alpha}-t_1^{1-4\alpha})^{1/2}<1$. Then
$\|u-v\|_{L^4(t_0,t_1,L_x^{\infty})}\le0$
which implies that $u(t)\equiv v(t)$ on 
$[t_0,t_1]$. This contradicts the assumption 
of $t_1$. Hence $u(t)=v(t)$ on $[T,\infty)$.

From (\ref{NLS}), we obtain
\begin{eqnarray}
u(t)=e^{-i(t-T)H_{q}}u(T)
-i\lambda\int_T^te^{-i\tau H_{q}}|u|^2u(\tau)d\tau.
\label{5.19}
\end{eqnarray}
Since $u(T)\in L_x^{2}(\rre)$, 
it follows from Lemma \ref{wp} that 
(\ref{5.19}) has a unique 
global solution in $C(\rre;L_x^2(\rre))\cap
L_{loc}^4(\rre;L_x^{\infty}(\rre))$. Therefore 
the solution $u$ of (\ref{NLS}) extends to all times. 
This completes the proof of Theorem \ref{main}. 
$\qquad\qed$

\section{Appendix} \label{sec:ape}

In this appendix, we give a proof of Lemma \ref{ly}. 
Since the case $y\ge0$ is trivial, 
we consider the case 
$y<0$ only. By definition, we have 
\begin{eqnarray*}
\lefteqn{{{\mathcal L}}_{+}[{{\mathcal F}}_{q}^{-1}[\phi]](y)}\\
&=&
{{\mathcal F}}_{q}^{-1}[\phi](y)
-qe^{qy}
\int_{y}^{0}e^{-qz}{{\mathcal F}}_{q}^{-1}[\phi](z)dz
-qe^{qy}
\int_{0}^{-y}e^{qz}{{\mathcal F}}_{q}^{-1}[\phi](z)dz
\end{eqnarray*}
Using the representation formula (\ref{J}) 
for ${{\mathcal F}}_{q}^{-1}$, 
we obtain 
\begin{eqnarray*}
\lefteqn{{{\mathcal L}}_{+}[{{\mathcal F}}_{q}^{-1}[\phi]](y)}\\
&=&{{\mathcal F}}^{-1}[\phi](y)
-qe^{-qy}\int_{-\infty}^{y}e^{qz}{{\mathcal F}}^{-1}[1_{+}\phi](z)dz\\
& &-qe^{-qy}\int^{\infty}_{-y}e^{-qz}
{{\mathcal F}}^{-1}[1_{-}\phi](z)dz\\
& &-qe^{qy}
\int_{y}^{0}e^{-qz}
{{\mathcal F}}^{-1}[\phi](z)dz
+q^{2}e^{qy}
\int_{y}^{0}e^{-2qz}
\left(\int_{-\infty}^{z}e^{qw}{{\mathcal F}}^{-1}[1_{+}\phi](w)dw
\right)dz\\
& &
+q^{2}e^{qy}
\int_{y}^{0}e^{-2qz}\left(\int^{\infty}_{-z}e^{-qw}
{{\mathcal F}}^{-1}[1_{-}\phi](w)dw\right)
dz\\
& &
-qe^{qy}
\int_{0}^{-y}e^{qz}
{{\mathcal F}}^{-1}[\phi](z)dz
+q^{2}e^{qy}
\int_{0}^{-y}e^{2qz}
\left(\int_{-\infty}^{-z}e^{qw}{{\mathcal F}}^{-1}[1_{+}\phi](w)dw
\right)dz\\
& &
+q^{2}e^{qy}
\int_{0}^{-y}e^{2qz}\left(\int^{\infty}_{z}e^{-qw}
{{\mathcal F}}^{-1}[1_{-}\phi](w)dw\right)
dz\\
&\equiv&G_{1}(x)+\cdots+G_{9}(x).
\end{eqnarray*}
Notice that for $y<0$,
\begin{eqnarray*}
\lefteqn{
e^{qy}\int_{y}^{0}e^{-2qz}
\left(
\int_{-\infty}^{z}e^{qw}f(w)dw\right)dz
=e^{qy}\int_{0}^{-y}e^{2qz}
\left(
\int_{-\infty}^{-z}e^{qw}f(w)dw\right)dz
}\\
&=&
\frac{1}{2q}e^{qy}\int_{y}^{0}e^{-qz}f(z)dz
+\frac{1}{2q}e^{-qy}\int_{-\infty}^{y}e^{qz}f(z)dz
-\frac{1}{2q}e^{qy}\int_{-\infty}^{0}e^{qz}f(z)dz,\\
\lefteqn{e^{qy}\int_{y}^{0}e^{-2qz}
\left(
\int_{-z}^{\infty}e^{-qw}f(w)dw\right)dz
=e^{qy}\int_{0}^{-y}e^{2qz}
\left(
\int_{z}^{\infty}e^{-qw}f(w)dw\right)dz
}\\
&=&
\frac{1}{2q}e^{qy}\int_{0}^{-y}e^{qz}f(z)dz
+\frac{1}{2q}e^{-qy}\int_{-y}^{\infty}e^{-qz}f(z)dz
-\frac{1}{2q}e^{qy}\int_{0}^{\infty}e^{-qz}f(z)dz.
\end{eqnarray*}
Hence we have
\begin{eqnarray*}
\lefteqn{G_{2}(x)+G_{5}(x)+G_{8}(x)}\\
&=&
qe^{qy}\int_{y}^{0}e^{-qz}{{\mathcal F}}^{-1}[1_{+}\phi](z)dz
-qe^{qy}\int_{-\infty}^{0}e^{qz}{{\mathcal F}}^{-1}[1_{+}\phi](z)dz
\\
&=&-G_{4}(x)
-qe^{qy}\int_{y}^{0}e^{-qz}{{\mathcal F}}^{-1}[1_{-}\phi](z)dz
-qe^{qy}\int_{-\infty}^{0}e^{qz}{{\mathcal F}}^{-1}[1_{+}\phi](z)dz,\\
\lefteqn{G_{3}(x)+G_{6}(x)+G_{9}(x)}\\
&=&
qe^{qy}\int_{0}^{-y}e^{qz}{{\mathcal F}}^{-1}[1_{-}\phi](z)dz
-qe^{qy}\int_{0}^{\infty}e^{-qz}{{\mathcal F}}^{-1}[1_{-}\phi](z)dz\\
&=&-G_{7}(x)
-qe^{qy}\int_{0}^{-y}e^{qz}{{\mathcal F}}^{-1}[1_{+}\phi](z)dz
-q\int_{0}^{\infty}e^{-qz}{{\mathcal F}}^{-1}[1_{-}\phi](z)dz.
\end{eqnarray*}
Collecting the above identities, we obtain (\ref{e1}). 
By a similar argument as above we also obtain (\ref{e2}),  
which completes the proof of Lemma \ref{ly}. 
$\qquad\qed$

\vskip2mm
\noindent
{\bf Acknowledgments.} The author is partially supported by MEXT, 
Grant-in-Aid for Young Scientists (A) 25707004 and 
the Sumitomo Foundation, Basic Science Research Projects No. 120043.




\begin{thebibliography}{30}

\bibitem{AG} 
Albeverio S., Gesztesy F., H\"{e}gh-Krohn R. and 
Holden, H.,
\textit{``Solvable Models in Quantum Mechanics''}. 
Springer-Verlag (1988).

\bibitem{B} 
Barab J. E.,
\textit{Nonexistence of asymptotically free solutions for a 
nonlinear Schr\"odinger equation}, 
J. Math. Phys., {\bf 25} (1984), 3270--3273.

\bibitem{C}
Cazenave T., \textit{``Semilinear Schr\"{o}dinger equations"}. 
Courant Lecture Notes in Mathematics, {\bf 10}.
American Mathematical Society (2003).

\bibitem{CGV}
Cuccagna S., Georgiev V. and Visciglia N., 
\textit{Decay and scattering of small solutions of pure power NLS 
in $\rre$ with $p>3$ and with a potential}. preprint 
(2012) To appear in  Comm. on Pure Appl. Math.

\bibitem{DH} Datchev K. and Holmer J., \textit{
Fast soliton scattering by attractive delta impurities}. 
Comm. Partial Differential Equations {\bf 34} (2009) 1074--1113.

\bibitem{DP} Deift P. and Park J., \textit{
Long-time asymptotics for solutions of the NLS equation 
with a delta potential and even initial data}. 
Int. Math. Res. Not. {\bf 2011} (2011), 5505--5624. 

\bibitem{FOO} Fukuizumi R., Ohta M. and Ozawa T., \textit{
Nonlinear Schr\"{o}dinger equation with a point defect}. 
Ann. Inst. H. Poincar\'{e} Anal. Non Linaire {\bf 25} 
(2008) 837--845.

\bibitem{GV3}
Ginibre J. and Velo G., \textit
{On a class of nonlinear 
Schr\"{o}dinger equations I. The Cauchy problem, 
general case. II. Scattering theory, general case.} 
J. Funct. Anal. {\bf 32} (1979), 1--71.

\bibitem{GHW}
Goodman R.H., Holmes P.J. and Weinstein M.I. 
\textit{Strong NLS soliton-defect interactions}. 
Phys. D {\bf 192} (2004) 215--248.


\bibitem{HN1} Hayashi N. and Naumkin P.I., 
\textit{Asymptotics for large time of
solutions to the nonlinear Schr\"{o}dinger and 
Hartree equations}. Amer. J. Math. 
{\bf 120} (1998), 369--389.

\bibitem{HN2} Hayashi N. and Naumkin P.I., 
\textit{Domain and range of the modified wave operator for 
Schr\"{o}dinger equations with a critical nonlinearity}. 
Comm. Math. Phys. 267 (2006), no. 2, 477--492.

\bibitem{HMZ1} Holmer J., Marzuola J. and Zworski M., 
\textit{Soliton splitting by external delta potentials}. 
J. Nonlinear Sci. {\bf 17} (2007) 349--367.

\bibitem{HMZ2} Holmer J., Marzuola J. and Zworski M., 
\textit{Fast soliton scattering by delta impurities}. 
Comm. Math. Phys. {\bf 274} (2007) 187--216.

\bibitem{HZ} Holmer J. and Zworski M., 
\textit{Slow soliton interaction with delta impurities}. 
J. Mod. Dyn. {\bf 1} (2007) 689--718.

\bibitem{KT}  Keel M. and Tao T., \textit{
Endpoint Strichartz estimates}. 
Amer. J. Math. {\bf 120} (1998), 955--980. 

\bibitem{LS}
Lin J. E. and Strauss W. A., \textit{
Decay and scattering of solutions of a 
nonlinear Schr\"{o}dinger equation}, 
J. Funct. Anal. {\bf 30} (1978), 245--263.

\bibitem{O}
Ozawa T., {\sl Long range scattering for 
nonlinear Schr\"{o}dinger equations in
one space dimension}. Comm. Math. Phys. 
{\bf 139} (1991), 479--493.

\bibitem{T} Tsutsumi Y., \textit{\ 
$L^2$-solutions for Nonlinear Schr\"odinger
equations and nonlinear groups}. 
Funkcialaj Ekvacioj. {\bf 30} (1987),
115--125.

\bibitem{TY}
Tsutsumi Y. and Yajima K.,  
{\sl The asymptotic behavior of nonlinear 
Schr\"{o}dinger equations}. 
Bull. Amer. Math. Soc. {\bf 11} (1984) 
186--188.

\bibitem{Y}
Yajima K.,  
\textit{Existence of solutions for Schr\"{o}dinger
evolution equations}. Comm. Math. Phys., {\bf 110} (1987), 415--426.
\end{thebibliography}
\end{document}